
\documentstyle{amsppt}
\baselineskip18pt
\magnification=\magstep1
\pagewidth{30pc}
\pageheight{45pc}
\hyphenation{co-deter-min-ant co-deter-min-ants pa-ra-met-rised
pre-print pro-pa-gat-ing pro-pa-gate
fel-low-ship Cox-et-er dis-trib-ut-ive}
\def\leaderfill{\leaders\hbox to 1em{\hss.\hss}\hfill}
\def\A{{\Cal A}}
\def\C{{\Cal C}}

\def\ldescent#1{{\Cal L (#1)}}
\def\rdescent#1{{\Cal R (#1)}}

\def\idest{i.e.,\ }

\def\ti{\widetilde}
\def\a{{a}}
\def\be{{b}}
\def\g{{c}}

\def\d{{\delta}}

\def\e{{\varepsilon}}

\def\bc{{\bold c}}

\def\b0{\text{\bf 0}}

\def\ra{{\ \longrightarrow \ }}

\def\zed{{\Bbb Z}}

\def\pd{\partial}
\def\boxit#1{\vbox{\hrule\hbox{\vrule \kern3pt
\vbox{\kern3pt\hbox{#1}\kern3pt}\kern3pt\vrule}\hrule}}
\def\rabbit{\vbox{\hbox{\kern0pt
\vbox{\kern0pt{\hbox{---}}\kern3.5pt}}}}

\def\tableau#1{
        \hbox {
                \hskip -10pt plus0pt minus0pt
                \raise\baselineskip\hbox{
                \offinterlineskip
                \hbox{#1}}
                \hskip0.25em
        }
}

\def\tabCol#1{
\hbox{\vtop{\hrule
\halign{\strut\vrule\hskip0.5em##\hskip0.5em\hfill\vrule\cr\lower0pt
\hbox\bgroup$#1$\egroup \cr}
\hrule
} } \hskip -10.5pt plus0pt minus0pt}

\def\CR{
        $\egroup\cr
        \noalign{\hrule}
        \lower0pt\hbox\bgroup$
}



\def\blank#1#2{
\hbox to #1{\hfill \vbox to #2{\vfill}}
}


\def\strut{\vrule height10pt depth5pt width0pt}

\topmatter
\title Acyclic heaps of pieces, II
\endtitle

\author R.M. Green \endauthor
\affil Department of Mathematics \\ University of Colorado \\
Campus Box 395 \\ Boulder, CO  80309-0395 \\ USA \\ {\it  E-mail:}
rmg\@euclid.colorado.edu \\
\endaffil

\abstract
We characterize and classify the ``regular classes of heaps'' introduced by
the author using ideas of Fan and of Stembridge.  The irreducible 
objects fall into five infinite families with one exceptional case.
\endabstract

\thanks
The author thanks Colorado State University for its hospitality during the
preparation of this paper.
\endthanks


\endtopmatter


\centerline{\bf To appear in the Glasgow Mathematical Journal}

\def\hker#1{\ker \pd_{#1}}

\head Introduction \endhead

A heap is an isomorphism class of labelled posets satisfying certain 
axioms.  Heaps have a wide variety of applications, as discussed by Viennot 
in \cite{{\bf 13}}.  

In \cite{{\bf 9}}, the author studied combinatorial properties
(called properties P1 and P2) which may or may not hold for a given heap;
property P1 is based on Fan's (algebraic) notion of cancellability 
\cite{{\bf 3}, \S4} and property P2 is related to Stembridge's definition of 
full commutativity \cite{{\bf 12}, \S1}.  When certain results of Fan from
\cite{{\bf 3}} and \cite{{\bf 4}, \S3} are reinterpreted in the context 
of heaps, we find that in certain heap monoids $H(P, \C)$ property P2 implies 
property P1.  In this case, we call $H(P, \C)$ a ``regular class of heaps''.
The basic combinatorial and linear properties of regular classes of heaps 
were developed in \cite{{\bf 9}, \S2}.  

The combinatorial properties of regular classes of heaps were used by Fan 
\cite{{\bf 3}, \S6} to investigate the representation theory of 
certain Hecke algebra quotients (also known as generalized Temperley--Lieb 
algebras).  Although slightly more subtle, the linear properties of regular 
classes of heaps are closely related to the combinatorial properties.
This theory was first developed in many important cases by Graham in his thesis
\cite{{\bf 5}} using a direct combinatorial argument, but our approach
seems to make the proofs more transparent as well as more general.  Graham
used these properties to obtain results on structure constants for the
Kazhdan--Lusztig bases of certain Hecke algebras \cite{{\bf 5}, \S9}.
In \cite{{\bf 10}}, J. Losonczy and the author used the same properties to
prove that in certain cases, the Kazhdan--Lusztig type bases of certain 
generalized Temperley--Lieb algebras are given by monomials in the generators.
The theory can also be applied to certain diagram calculi for
these algebras.  For further applications and more details, the reader is
referred to \cite{{\bf 9}, \S4.1}.

In the light of these properties, it is desirable to obtain a better
understanding of regular classes of heaps.  The first main result of this 
paper (Theorem 1.5.1, proved in \S2) gives some equivalent characterizations 
of the property of being regular; these involve linear algebra and certain 
associative algebras as well as combinatorial properties.  The second main 
result (Theorem 1.5.2, proved in \S3) solves the
problem posed in \cite{{\bf 9}, Problem 4.3.1} and gives a complete 
classification of regular classes of heaps, assuming the corresponding
set of pieces is finite.   These are classified by their concurrency
graphs, and the irreducible objects fall into five infinite families together 
with one exceptional case.

\head 1. Preliminaries \endhead

We begin with some preliminary material that is necessary for the statement
of the main results.  Our approach follows \cite{{\bf 13}} and \cite{{\bf 9}}.

\subhead 1.1 Heaps \endsubhead

We start by recalling the basic definitions following the conventions of
\cite{{\bf 9}}.  These differ slightly from those of \cite{{\bf 13}}; see \cite{{\bf 9},
\S1.1} for details.

\definition{Definition 1.1.1}
Let $P$ be a set equipped with a symmetric and reflexive binary relation
$\C$.  The elements of $P$ are called {\it pieces}, and the relation
$\C$ is called the {\it concurrency relation}.

A {\it labelled heap} with pieces in $P$ is a triple $(E, \leq, \e)$ 
where $(E, \leq)$ is a finite (possibly empty)
partially ordered set with order relation denoted
by $\leq$ and $\e$ is a map $\e : E \ra P$ satisfying the following two
axioms. 

\item{1.}{For every $\a, \be \in E$ such that $\e(\a) \ \C \ \e(\be)$, 
$\a$ and $\be$ are comparable in the order $\leq$.}

\item{2.}{The order relation $\leq$ is the transitive closure of the
relation $\leq_\C$ such that for all $\a, \be \in E$, $\a \ \leq_\C \ \be$ 
if and only if both $\a \leq \be$ and $\e(\a) \ \C \ \e(\be)$.}
\enddefinition

The terms {\it minimal} and 
{\it maximal} applied to the elements of the labelled 
heap refer to minimality (respectively, maximality) with respect to $\leq$.

\example{Example 1.1.2}
Let $P = \{1, 2, 3\}$ and, for $x, y \in P$, define $a \ \C \ b$ 
if and only if 
$|x - y| \leq 1$.  Let $E = \{a, b, c, d, e\}$ partially ordered by
extension of the (covering) relations $a \leq c$, $b \leq c$, $c \leq d$,
$c \leq e$.  Define the map $\e$ by the conditions $\e(a) = \e(d) = 1$,
$\e(c) = 2$ and $\e(b) = \e(e) = 3$.  Then $(E, \leq, \e)$ can easily be
checked to satisfy the axioms of Definition 1.1.1 and it is a labelled heap.
The minimal elements are $a$ and $b$, and the maximal elements are $d$ and $e$.
\endexample

\definition{Definition 1.1.3}
Let $(E, \leq, \e)$ and $(E', \leq', \e')$ be two labelled 
heaps with pieces in $P$ and with the same concurrency relation, $\C$.  
An isomorphism $\phi : E \ra E'$ of posets is said to be an 
{\it isomorphism of labelled posets} if $\e = \e' \circ \phi$.

A {\it heap} of pieces in $P$ with concurrency relation $\C$ is a labelled
heap (Definition 1.1.1) defined up to labelled poset isomorphism.
The set of such heaps is denoted by $H(P, \C)$.  We denote the heap 
corresponding to the labelled heap $(E, \leq, \e)$ by $[E, \leq, \e]$.
\enddefinition

We will sometimes abuse language and speak of the underlying set of a heap,
when what is meant is the underlying set of one of its representatives.

\definition{Definition 1.1.4}
Let $(E, \leq, \e)$ be a labelled heap with pieces in $P$ and $F$ a 
subset of $E$.  
Let $\e'$ be the restriction of $\e$ to $F$.  Let ${\Cal R}$ be the relation
defined on $F$ by $\a \ {\Cal R} \ \be$ if and only if $\a \leq \be$ and
$\e(\a) \ \C \ \e(\be)$.  Let $\leq_F$ be the transitive closure of ${\Cal R}$.
Then $(F, \leq_F, \e')$ is a labelled heap with pieces in $P$.  The heap 
$[F, \leq_F, \e']$ is called a {\it subheap} of $[E, \leq, \e]$.
If $v \in E$, we let $E(v) = [E(v), \leq_{E(v)}, \e']$ be the subheap of $E$ 
obtained by defining $E(v) = E \backslash \{v\}$.
\enddefinition

We will often implicitly use the fact that a subheap is determined by its
set of vertices and the heap it comes from.

\definition{Definition 1.1.5}
The {\it concurrency graph} associated to the class of heaps $H(P, \C)$ is
the graph whose vertices are the elements of $P$ and for which there is an
edge from $v \in P$ to $w \in P$ if and only if $v \ne w$ and $v \ \C \ w$.
\enddefinition

\definition{Definition 1.1.6}
Let $E = [E, \leq_E, \e]$ and $F = [F, \leq_F, \e']$ be two heaps in 
$H(P, \C)$.
We define the heap $G = [G, \leq_G, \e''] = E \circ F$ of $H(P, \C)$ 
as follows.

\item{1.}{The underlying set $G$ is the disjoint union of $E$ and $F$.}
\item{2.}{The labelling map $\e''$ is the unique map $\e'' : G \ra P$ whose
restriction to $E$ (respectively, $F$) is $\e$ (respectively, $\e'$).}
\item{3.}{The order relation $\leq_G$ is the transitive closure of the
relation ${\Cal R}$ on $G$, where $\a \ {\Cal R} \ \be$ if
and only if one of the following three conditions holds:
\item{(i)}{$\a, \be \in E$ and $\a \leq_E \be$;}
\item{(ii)}{$\a, \be \in F$ and $\a \leq_F \be$;}
\item{(iii)}{$\a \in E, \ \be \in F$ and $\e(\a) \ \C \ \e'(\be)$.}
}
\enddefinition

\remark{Remark 1.1.7}
Definition 1.1.6 can easily be shown to be sound (see \cite{{\bf 13}, \S2}).
It is immediate from the construction that $E$ and $F$ are subheaps of 
$E \circ F$.

As in \cite{{\bf 13}}, we will write $\a \circ E$ and $E \circ \a$ for
$\{\a\} \circ E$ and $E \circ \{ \a \}$, respectively.  Note that $\a \circ E$
and $\be \circ E$ are equal as heaps if $\e(\a) = \e(\be)$.
If $I$ is a finite set, we may also write
$\prod_{i \in I} \a_i$ for the product of the heaps $\{\a_i\}$ in the case
where the singleton heaps commute pairwise.
\endremark

\definition{Definition 1.1.8}
A {\it trivial heap} is a heap $[E, \leq, \e]$ for which the order relation
$\leq$ is trivial.
\enddefinition

\subhead 1.2 Convex chains and regular classes \endsubhead

In \S1.2, 
we recall the definitions of properties P1 and P2 
for heaps from \cite{{\bf 9}, \S2.2}, and 
we recall some terminology relating to convex chains from 
\cite{{\bf 9}, \S2.3}.

\definition{Definition 1.2.1}
A {\it convex chain} in a heap $E = [E, \leq, \e]$ is a chain  $$
\bc = (x_1, x_2, \ldots, x_t) : x_1 < x_2 < \cdots < x_t
$$ of vertices in $E$ such that whenever $x_i < y < x_j$ for some
$y$, the vertex $y$ is an element of the chain.  A convex chain is said to
be {\it balanced} if $\e(x_1) = \e(x_t)$.  If $\bc$ is a balanced 
convex chain, we define the heap $E/\bc$ to be the subheap of $E$ 
obtained by omitting the vertices $x_2, x_3, \ldots x_t$.  We call the 
heap $E/\bc$ the {\it contraction of $E$ along $\bc$}, and the
number $t$ is called the {\it length} of the chain.
\enddefinition

We can improve on the notation of \cite{{\bf 9}} by using the following 
notation (based on \cite{{\bf 3}, \S2}) for maximal and minimal elements of a heap.

\definition{Definition 1.2.2}
Let $E = [E, \leq, \e]$ be a heap.  We define $\ldescent{E}$ to be the set of 
vertices minimal in $E$, and $\rdescent{E}$ to be the set of vertices 
maximal in $E$.
\enddefinition

\example{Example 1.2.3}
If $E$ is the heap arising from Example 1.1.2, we have $\ldescent{E} 
= \{a, b\}$ and $\rdescent{E} = \{c, d\}$.
\endexample

\definition{Definition 1.2.4 (Property P1)}
Let $E = [E, \leq, \e] \in H(P, \C)$ be a heap.
We write $E(\a) \prec^+ E$ (respectively, $E(\a) \prec^- E$) 
if $\a \in \rdescent{E}$ (respectively, $\a \in \ldescent{E}$) and there
exists $\be \in \rdescent{E(\a)} \backslash \rdescent{E}$ (respectively, 
$\be \in \ldescent{E(\a)} \backslash \ldescent{E}$) with
$\e(\be) \ne \e(\a)$.  We write $E(\a) \prec E$ if either
$E(\a) \prec^+ E$ or $E(\a) \prec^- E$.

If there is a (possibly trivial) sequence $
E_1 \prec E_2 \prec \cdots \prec E
$ of heaps in $H(P, \C)$ where $E_1$ is a trivial heap, we say that 
the heap $E$ is {\it dismantlable} or that $E$ has {property P1}.
\enddefinition

\remark{Remark 1.2.5}
Property P1 is closely related to Fan's notion of
left and right cancellability \cite{{\bf 3}, Definition 4.2.4}.  
As in \cite{{\bf 9}}, we avoid
the term ``cancellability'' in this paper because of possible confusion 
with the use of this term in the theory of monoids.
\endremark

\definition{Definition 1.2.6 (Property P2)}
We say a heap $E = [E, \leq, \e] \in H(P, \C)$
has property P2 if it contains no balanced convex chains of length 2 or 3.
\enddefinition

\remark{Remark 1.2.7}
Property P2 is modelled on Stembridge's characterization of full
commutativity \cite{{\bf 12}, Proposition 2.3} in the case of a 
simply laced Coxeter group.
\endremark

\subhead 1.3 Acyclic heaps \endsubhead

We recall the definition of the map $\pd$ from \cite{{\bf 9}, \S1.2}, to which
the reader is referred for further elaboration and examples.
Throughout \S1.3, we let $[E, \leq, \e]$ be a heap
in the set $H(P, \C)$ with pieces in $P$ and concurrency relation $\C$.  
We also fix a field, $k$.

\definition{Definition 1.3.1}
Let $V_0$ be the set of elements of $[E, \leq, \e]$, \idest the set of 
elements of (a representative of) the underlying poset, $E$.
We call the elements of $V_0$ {\it vertices} and  denote their $k$-span 
by $C_0$.

Let $V_1$ be the set of all pairs $(x, y) \in E \times E$ with $x < y$ and 
$\e(x) = \e(y)$ such that there is no element $z$ for which we have
both $\e(x) = \e(z) = \e(y)$ and $x < z < y$.  
We call the elements of $V_1$ {\it edges} and denote their $k$-span by $C_1$.

The $k$-linear map $\pd = \pd_E : C_1 \ra C_0$ is defined by its
effect on the edges as follows: $$
\pd : (x, y) \mapsto \sum_{{x < w < y} \atop {\e(w) \ \C \ \e(x)}} w
.$$
\enddefinition

\definition{Definition 1.3.2}
Let $E = [E, \leq_E, \e]$ be a heap in $H(P, \C)$ and let $k$ be a field.  
We say $E$ is {\it acyclic} if $\hker{E} = 0$.  We say $E$ is {\it strongly
acyclic} if $E$ is acyclic and $E(v)$ is acyclic for all $v \in E$.
\enddefinition

Some of the main results of \cite{{\bf 9}} may be summarized in the following
theorem.

\proclaim{Theorem 1.3.3 \cite{{\bf 9}}}
\item{\rm (a)}{Any strongly ayclic heap has property P2.}
\item{\rm (b)}{Any heap with property P1 is acyclic.}
\item{\rm (c)}
{Let $H(P, \C)$ be a regular class of heaps, and let $E = [E, \leq, \e]$ 
be a heap of $H(P, \C)$.  Then:
\item{\rm (i)}{$E$ has property P2 if and only if it is strongly acyclic;}
\item{\rm (ii)}{$E$ has property P1 if and only if it is acyclic.}
}
\endproclaim

\demo{Proof}
For (a), see \cite{{\bf 9}, Proposition 2.2.7}.  For (b), see \cite{{\bf 9}, 
Proposition 2.2.3}.  For (c), see \cite{{\bf 9}, Theorem 2.4.2} for (i) and
\cite{{\bf 9}, Theorem 2.4.4} for (ii).
\qed\enddemo

\remark{Remark 1.3.4}
Our characterization of regular classes of heaps in Theorem 1.5.1 will 
show that
the criterion in Theorem 1.3.3 (c)(i) is equivalent to $H(P, \C)$ being
regular.  The same is not true of Theorem 1.3.3 (c)(ii), although 
we do not prove this here.
\endremark

\subhead 1.4 Algebras arising from heaps \endsubhead

It is sometimes convenient to phrase arguments about the combinatorics or
linear properties of heaps in an algebraic way.  For this, we recall
the definition of certain monoids and algebras discussed in 
\cite{{\bf 9}, \S3.1}.

\definition{Definition 1.4.1}
A class of heaps $H(P, \C)$ has a natural monoid structure with composition
given by the map $\circ$ of Definition 1.1.6.  We call this
monoid the {\it heap monoid}.
\enddefinition

Another way to approach the heap monoid is by considering the commutation
monoids of Cartier and Foata \cite{{\bf 1}}, which are defined as follows.

\definition{Definition 1.4.2}
Let $A$ be a set and let $A^*$ be the free monoid generated by $A$.  Let $C$
be a symmetric and antireflexive relation on $A$.  The 
{\it commutation monoid} $\text{\rm Co} (A, C)$ is the quotient of the free
monoid $A^*$ by the congruence $\equiv_C$ generated by the commutation
relations: $$
ab \equiv_C ba \text{ for all } a, b \in A \text{ with } a \  C \ b
.$$
\enddefinition

The following result, proved in \cite{{\bf 13}, Proposition 3.4}, 
shows that the heap monoid is
naturally isomorphic to a commutation monoid.

\proclaim{Proposition 1.4.3}
Let $H(P, \C)$ be a class of heaps and let $C$ be the 
complementary relation of $\C$.  Let $E = [E, \leq, \e]$ be a heap of 
$H(P, \C)$.  The map from $H(P, \C)$ to $\text{\rm Co}(P, C)$ that sends the 
heap $$
\a_1 \circ \a_2 \circ \cdots \circ \a_r \mapsto 
\e(\a_1) \e(\a_2) \cdots \e(a_r)
$$ is an isomorphism of monoids. \qed
\endproclaim

It is convenient to study certain quotients of heap monoid algebras.  The
algebras below and their bases appear in the work of Fan \cite{{\bf 2}} and
Graham \cite{{\bf 5}}.

\definition{Definition 1.4.4}
Maintain the above notation, so that $C$ is the complementary relation of
$\C$.  Let $\A$ be the ring of
Laurent polynomials $\zed[v, v^{-1}]$, let $\d := v + v^{-1}$, 
and let $\A \text{\rm Co}(P, C)$ be the monoid
algebra of $\text{\rm Co}(P, C)$ over $\A$.  We define the {\it generalized 
Temperley--Lieb algebra} $TL(P, \C)$ to be the $\A$-algebra 
obtained by quotienting $\A \text{\rm Co}(P, C)$ by the relations $$\eqalign{
s s &= \d s,\cr
s t s &= s \text{ if } s \ne t \text{ and } s \ \C \ t,\cr
}$$ where $s, t \in P$.
\enddefinition

\proclaim{Lemma 1.4.5}
The isomorphism of Proposition 1.4.3 induces an isomorphism between the 
algebra $TL(P, \C)$ and the quotient of $\A H(P, \C)$ by 
the relations $$\eqalign{
E &= \d E/\bc \text{ if } \bc \text{ is a balanced convex chain 
of length 2},\cr
E &= E/\bc \text{ if } \bc \text{ is a balanced convex chain } x < y < z
\text{ with } \e(x) \ne \e(y).\cr
}$$
\endproclaim

\demo{Proof}
See \cite{{\bf 9}, Lemma 3.1.6}.
\qed\enddemo

\proclaim{Proposition 1.4.6}
The quotient of $\A H(P, \C)$ described in Lemma 1.4.5 has as a free
$\A$-basis the images of those heaps in $H(P, \C)$ with property P2.
\endproclaim

\demo{Proof}
See \cite{{\bf 9}, Proposition 3.2.2}.
\qed\enddemo

\definition{Definition 1.4.7}
The basis of $TL(P, \C)$ corresponding to the basis of $\A H(P, \C)$ given
in Proposition 1.4.6 under the isomorphism of Proposition 1.4.3 is called the
{\it monomial basis} of $TL(P, \C)$.  (It consists of certain monomials in 
the set $P$.)
\enddefinition

\subhead 1.5 Statement of main results \endsubhead

The main two results of this paper concern the characterization and 
classification of regular classes of heaps.  \S2 will be devoted to Theorem
1.5.1, and \S3 will be devoted to Theorem 1.5.2.

\proclaim{Theorem 1.5.1}
Let $H(P, \C)$ be the class of heaps with pieces in $P$ and concurrency 
relation $\C$.  The following are equivalent:
\item{\rm (i)}{every heap in $H(P, \C)$ with property P2 has property P1
(\idest $H(P, \C)$ is regular);}
\item{\rm (ii)}{every heap in $H(P, \C)$ with property P2 is strongly
acyclic;}
\item{\rm (ii$'$)}{property P2 and the property of being strongly acyclic are 
equivalent for heaps in $H(P, \C)$;}
\item{\rm (iii)}{if $E$ is a heap with property P2 then, for any vertex $\a$
of $E$, $E(\a)$ is equal as an element of $TL(P, \C)$ to a heap with property
P2;}
\item{\rm (iii$'$)}{whenever $p_1 p_2 \cdots p_r$ is a monomial basis 
element of $TL(P, \C)$, the monomial obtained by deleting one of the 
$p_i$ is equal (as an element of $TL(P, \C)$) to a monomial basis element.}
\endproclaim

\topinsert
\topcaption{Figure 1} Connected graphs associated with regular classes of heaps
\endcaption
\centerline{
\hbox to 3.319in{
\vbox to 4.638in{\vfill
        \includegraphics{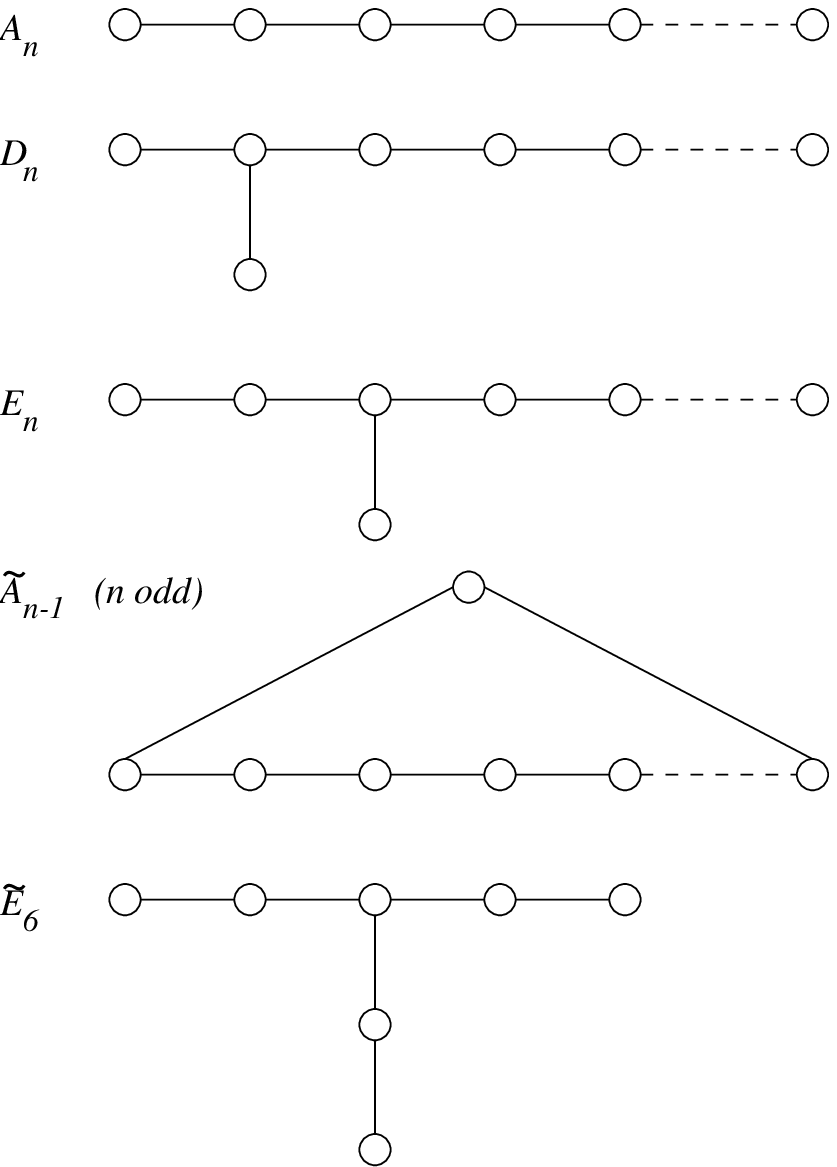}
}
\hfill}
}
\endinsert

The classification of regular classes of heaps is in terms of their 
concurrency graphs (see Definition 1.1.5).  This involves the graphs in
Figure 1.  For types $A_n$, $D_n$ and $E_n$, $n$ is the number of nodes in
the graph, and we assume $n \geq 4$ for type $D_n$ and $n \geq 6$ for type
$E_n$ to avoid repetition.  (Note that we require graphs $E_n$ for arbitrarily
high $n$.)  The graph of type $\ti{A}_{n-1}$ has $n$ nodes where $n$ is odd
and $n \geq 3$.  The graph of type $\ti{E}_6$ has exactly 7 nodes.

\proclaim{Theorem 1.5.2}
Let $H(P, \C)$ be the class of heaps with pieces in $P$ and concurrency 
relation $\C$.  Suppose also that $P$ is a finite set, and let $\Gamma$ be
the concurrency graph associated to $H(P, \C)$.  Then $H(P, \C)$ is a regular
class of heaps if and only if each connected component of $\Gamma$ is
a complete graph $K_n$ or appears
in the list depicted in Figure 1: type $A_n$ ($n \geq 1$), type $D_n$ 
($n \geq 4$), type $E_n$ ($n \geq 6$), type $\ti{A}_{n-1}$ ($n \geq 3$ and
$n$ odd) or type $\ti{E}_6$.
\endproclaim

\remark{Remark 1.5.3}
We assume that $P$ is finite above to avoid cardinality considerations.
\endremark

\head 2. Characterization of regular classes of heaps \endhead

The aim of \S2 is to prove Theorem 1.5.1.  In order to do so we shall need
to recall and develop some additional combinatorial properties of heaps.

\subhead 2.1 Factorization of heaps \endsubhead

Heaps possess the following unique factorization property.

\proclaim{Proposition 2.1.1}
Any heap $E \in H(P, \C)$ can be written uniquely as a product of trivial
heaps $$
E = T_1 \circ T_2 \circ \cdots \circ T_p
$$ such that for each $1 \leq j < p$ and for each $\be \in T_{j+1}$, 
there exists $\a \in T_j$ with $\a < \be$.
\endproclaim

\demo{Note}
Note that, in the above situation, it is possible for $\e(\a) = \e(\be)$.
\enddemo

\demo{Proof}
See \cite{{\bf 13}, Lemma 2.9}.
\qed\enddemo

\example{Example 2.1.2}
The factorization of the heap $E$ arising from Example 1.1.2 is of the form 
$T_1 \circ T_2 \circ T_3$ whose underlying sets are given by $T_1 = \{a, b\}$,
$T_2 = \{c\}$ and $T_3 = \{d, e\}$.
\endexample

\remark{Remark 2.1.3}
In the unique factorization of the heap $E$ given by Proposition 2.1.1, 
the heap $T_1$ is the subheap of $E$ consisting of the vertices $\ldescent{E}$.
(See the proof of \cite{{\bf 13}, Lemma 2.9} for justification, and Definition
1.2.2 for the notation.)  It is not always the case that
the rightmost subheap in the factorization of $E$ coincides with $\rdescent{E}$
(although it does in Example 1.2.3).
\endremark

\definition{Definition 2.1.4}
Let $E = [E, \leq, \e]$ be a heap of $H(P, \C)$ and let $$
E = T_1 \circ T_2 \circ \cdots \circ T_p
$$ be its unique factorization as in Proposition 2.1.1.  We define
$E^* = [E, \geq, \e]$, the {\it opposite heap} of $E$, to be the heap $$
T_p \circ T_{p-1} \circ \cdots \circ T_1
.$$  If $\a$ is a vertex of $E$ situated in the factor $T_i$ then we denote
the corresponding element of $E^*$ in the factor $T_i$ by $\a^*$.  (Note
that $\e(\a^*) = \e(\a)$.)

We define the 
{\it double}, $\Delta(E)$, of $E$ to be the heap of $H(P, \C)$ given by $$
\Delta(E) = T_p \circ T_{p-1} \circ \cdots \circ T_2 \circ T_1
\circ T_2 \circ \cdots \circ T_{p-1} \circ T_p
.$$  Note that $\Delta(E)$ contains both $E$ and $E^*$ as subheaps.
\enddefinition

\remark{Remark 2.1.5}
The expression given above for $\Delta(E)$ in terms of trivial heaps may or
may not be the unique factorization of $\Delta(E)$ given by Proposition 2.1.1.
\endremark

\example{Example 2.1.6}
Consider the heap of Example 1.1.2.
In this case, the opposite heap $E^*$ is given by $$
(d^* \circ e^*) \circ (c^*) \circ (a^* \circ b^*)
,$$ where the parentheses enclose the factors in the unique factorization
of $E^*$.  The double of $E$ is given by $$
\Delta(E) =  
(d^* \circ e^*) \circ (c^*) \circ
(a \circ b) \circ (c) \circ (d \circ e) 
,$$ which happens to be the factorization of $\Delta(E)$ given
by Proposition 2.1.1.
\endexample

\subhead 2.2 More on the relation $\prec$ \endsubhead

For the main results of \S2, we will need to examine more closely 
the relation $\prec$ that appears in the definition of
property P1.  Interesting examples of these results will be found in \S3
when we classify the regular classes of heaps.

The next result is used in \cite{{\bf 3}, \S4}.

\proclaim{Lemma 2.2.1}
Let $E = [E, \leq, \e]$ be a heap with unique factorization (as in Proposition
2.1.1) $$
E = T_1 \circ T_2 \circ \cdots \circ T_p
$$ where $p > 1$, and suppose there is no subheap $F$ of $E$ with 
$F \prec^- E$.  Then for each vertex $\be \in T_2$, either 
\item{\rm (i)}{there is a vertex
$\a \in T_1$ such that $\e(\a) = \e(\be)$, or}
\item{\rm (ii)}{there are at least
two vertices $\a_1, \a_2 \in T_1$ such that
$\a_i < \be$ and such that $\e(\a_1), \e(\a_2), \e(\be)$ are distinct.}
\endproclaim

\demo{Proof}
Suppose that (i) does not hold.  (Recall that the vertices in $T_1$ are
precisely the minimal vertices of $E$.)

The existence of one vertex, $\a_1$, satisfying the hypotheses is 
guaranteed by Proposition 2.1.1.  Since (i) does not hold, $\e(\a_1) \ne
\e(\be)$, which means that if there is no vertex $\a_2$ with $\e(\a_2) \ne
\e(\be)$ that is distinct from $\a_1 \in T_1$, we must
have $\be \in \ldescent{E(\a_1)} \backslash \ldescent{E}$ and $E(\a_1)
\prec^- E$, a contradiction.  Finally, $\e(\a_1) \ne \e(\a_2)$, because
there cannot be a chain $\a_1 < \a_2 < \be$ unless $\be \in T_i$ for
$i \geq 3$.
\qed\enddemo

\proclaim{Lemma 2.2.2}
Maintain the notation of Definition 1.2.4.  Let $E = [E, \leq, \e]$ be a
heap of $H(P, \C)$ containing no balanced convex chains of length 2.  
Suppose that there is no subheap $F$ of $E$ with $F \prec^+ E$.
\item{\rm (i)}{If $E(\a) \prec^- E$, then $
\Delta(E) = \d^m \Delta(E(\a))
$ as elements of $TL(P, \C)$, for some nonnegative integer $m$.
If, in addition, $E(\a)$ is trivial, we have $m \geq 1$.}
\item{\rm (ii)}{If $\a$ is a minimal element of $E$ (in particular,
if $E(\a) \prec^- E$) then
there is no subheap $F'$ of $E(\a)$ with $F' \prec^+ E(\a)$.}
\endproclaim

\demo{Proof}
We first prove (ii).  Definition 1.2.4 implies that if $F' \prec^+ E(\a)$
then $\a \circ F' \prec^+ E$: by hypothesis, $\a$
does not interfere with the properties of maximal elements required by the
relation $\prec^+$, proving (ii).

Suppose now that $E(\a) \prec^- E$ as in the statement of (i), so that
there is a vertex $\be \in \ldescent{E(\a)} 
\backslash \ldescent{E}$ with $\e(\be) \ne \e(\a)$.  We denote the (possibly
empty) set $\ldescent{E(\a)} \backslash (\ldescent{E} \cup \{ \be \})$ by $M$.
Let $T_1 \circ T_2 \circ \cdots \circ T_p$ be the unique factorization
of the heap $E(\a)$.

Now $\be^* \circ \a \circ \be$ gives a balanced convex chain, 
$\bc$, of length 3 in the double $\Delta(E)$ of $E$.  Since
$\e(\be^*) = \e(\be) \ne \e(\a)$, Lemma 1.4.5 shows that $\Delta(E) = 
\Delta(E)/\bc$ in $TL(P, \C)$.  The heap $\Delta(E)/\bc$ is closely related
to $\Delta(E(\a))$; the only difference is that each element $\mu$ of $T_1$ in
$\Delta(E(\a))$ corresponding to an element of the set $M$ is replaced in
$\Delta(E)/\bc$ by two elements, $\mu^* \circ \mu$.  These correspond to
balanced convex chains of length 2 and we may apply Lemma 1.4.5 again to show
that $\Delta(E) = \d^m \Delta(E(\a))$, where $m = |M|$.

It remains to consider the case where $E(\a)$ is trivial; this implies that $E$
is not trivial because $E(\a) \prec^- E$.  In this case,
$E = \a \circ E(\a)$.  Let $T'_1 \circ T'_2$ be the unique factorization
of $E$; then the vertex $\a$ lies in the factor $T'_1$.
Let $T''_1 \circ T''_2$ be the unique factorization
of the opposite heap, $E^*$ of $E$; the factor $T''_2$ contains only the vertex
$\a^*$.  Since (by hypothesis) there is no subheap $F$ of $E$
with $F \prec^+ E$, it follows that there is no subheap $F$ of $E^*$ with
$F \prec^- E^*$.  Lemma 2.2.1 applies to the heap $E^*$, but case (i) cannot
occur because $E$ (and therefore $E^*$) contains no balanced 
convex chains of length 2.  There are therefore at least two vertices 
$\be_i^*$ in $E^*$ such that $\be_i^* < \a^*$, meaning that $\be_i >
\a$ for each such vertex.  The vertices $\be_i \in E$ all lie in the 
factor $T'_2$, which is a trivial subheap.  Let $M'$ be the set of all
such vertices $\be_i$.

Iterated applications of Lemma 1.4.5 show that $$
\left( \prod_{\be \in M'} \be^* \right) \circ \a \circ 
\left( \prod_{\be \in M'} \be \right) = \d^{|M'| - 1} 
\left( \prod_{\be \in M'} \be \right)
$$ as elements of $TL(P, \C)$.  (See Remark 1.1.7 for the notation.)  Since
the elements of $E$ that do not appear in the above product (a) commute with 
the elements that do appear and (b) appear in the factor
$T'_1$, we have $$
\Delta(E) = \d^{|M'| - 1} \Delta(E(\a))
,$$ and the conclusion of (i) follows from the fact that $|M'| \geq 2$.
\qed\enddemo

\proclaim{Lemma 2.2.3}
Maintain the above notation.  Suppose that $E = [E, \leq, \e]$ is a heap of
$H(P, \C)$ with property P1 and that there is no subheap $F$ of $E$ with
$F \prec^+ E$.  Then:
\item{\rm (i)}{there is a sequence $$
E_1 \prec^- E_2 \prec^- \cdots \prec E_{r-1} \prec^- E_r = E
$$ of heaps with $E_1$ trivial, and for each $E_i$ there is no subheap $F_i$
of $E_i$ with $F_i \prec^+ E_i$;}
\item{\rm (ii)}{$\Delta(E) = \d^m E_1$ as elements of $TL(P, \C)$ for some
nonnegative integer $m$, and $m > 0$ if $E$ is nontrivial;}
\item{\rm (iii)}{if $\a$ is a minimal element of $E$ then there is a
sequence like that of (i) but with final term $E_r = E(\a)$; 
in particular, $E(\a)$
has property P1.}
\endproclaim

\demo{Proof}
Part (i) is a consequence of the definition of property P1 and Lemma 
2.2.2 (ii), whose hypotheses are satisfied by Theorem 1.3.3 (b).

For part (ii), note that $E_1 = \Delta(E_1)$ because $E_1$ is trivial.  The
claim follows from (i) and repeated applications of Lemma 2.2.2 (i); if $E$ is
nontrivial, we have $m > 0$ by considering $E_1 \prec^- E_2$.

We now turn to (iii).  The steps $E_{i-1} \prec^- E_i$ in the sequence in 
(i) correspond to the removal of a sequence of (distinct) vertices $\a_i$, 
each minimal in $E_i$.  According to the definition of $\prec^-$,
the removal of each vertex $\a_i$ exposes a new minimal vertex $\be_i$ that is
not minimal in $E_i$.  Since $\a$ is minimal in $E$ we cannot have 
$\a_i < \a$ for any of the elements $\a_i$.  

There are now two cases to consider.  In the first case, $\a$ is not equal
to any of the $\a_i$.  In this case, we can adapt the sequence in (i) by
replacing each $E_i$ by $E_i(\a)$, and the required properties hold.  In
the second case, $\a$ is equal to one of the $\a_i$.  In this case,
the sequence $$
E_1 \prec^- E_2 \prec^- \cdots \prec^- E_{i-1} = E_i(\a) 
\prec^- E_{i+1}(\a) \prec^- 
\cdots \prec^- E_r(\a) = E(\a)
$$ has the required properties.
\qed\enddemo

\subhead 2.3 More on doubles of heaps \endsubhead

\proclaim{Lemma 2.3.1}
If $E = [E, \leq, \e]$ is a heap of $H(P, \C)$ with property P2
but not property P1, then $\Delta(E)$ has property P2.
\endproclaim

\demo{Proof}
It is enough to show that if $x < z$ are two vertices of $\Delta(E)$ with
$\e(x) = \e(z)$, then
there exist at least two distinct vertices $y_1$ and $y_2$ with $x < y_1 < z$,
$x < y_2 < z$ and $\e(y_i) \ne \e(x)$ for $i \in \{1, 2\}$.  
We may assume that there is no element $w$ with $x < w < z$ and 
$\e(x) = \e(w) = \e(z)$, or $E$ would have a balanced convex chain of
length $2$.

Let us write $$
\Delta(E) = T_p \circ T_{p-1} \circ \cdots \circ T_2 \circ T_1
\circ T_2 \circ \cdots \circ T_{p-1} \circ T_p
.$$  Note that $p > 1$ because $E$ cannot be trivial.

If $x$ and $z$ both come from factors in or to the right of the factor
$T_1$, then we are done because $E$ has property P2.  Similarly, if $x$ and
$z$ both come from factors in or to the left of the factor $T_1$, we are
done because $E^*$ inherits property P2 from $E$.  We may therefore assume
that neither $x$ nor $z$ comes from $T_1$, that $x$ comes from a factor to
the left of $T_1$, and that $z$ comes from a factor to the right of $T_1$.
The assumptions made above regarding the element $w$ imply that $x = z^*$.

Suppose that $z$ comes from the factor $T_i$ with $i > 2$.  This means that $z$
is not minimal in $E$.  By Proposition 2.1.1, there is an element $y$ in the
factor $T_{i-1}$ with $y < z$ and $\e(y) \ne \e(z)$ (since $E$ has 
property P2), and we may take $y_1 = y$ and
$y_2 = y^*$ (note that $y^* \ne y$ because $i > 2$).

The other case is when $z$ comes from the factor $T_2$.  In this case
we apply Lemma 2.2.1 with $\be = z$.  Case (i) cannot hold, because it would
contradict the assumption that $E$ has property P2.  We then take $y_i = \a_i$
for $i \in \{1, 2\}$.
\qed\enddemo

\demo{Proof of Theorem 1.5.1}
The equivalence of (ii) and (ii$'$) of Theorem 1.5.1 follows from Theorem 
1.3.3 (a).  The equivalence of (iii) and (iii$'$) follows from Definition
1.4.7.

The implication (i) $\Rightarrow$ (ii) comes from \cite{{\bf 9},
Theorem 2.4.2} (see Theorem 1.3.3 (c) (i)) for the statement).

Next consider the implication (ii) $\Rightarrow$ (iii).  By \cite{{\bf 9},
Theorem 3.2.3}, we see that $E(\a) = \d^m G$ in $TL(P, \C)$, 
for some heap $G$ with property
P2 and some nonnegative integer $m$.  Furthermore, that result shows
that $\dim \hker{E(\a)} = m + \dim \hker{G}$.  Since, by (ii), $G$ is
(strongly) acyclic, we have $\dim \hker{E(\a)} = m$.  Since $E$ has property
P2, it is strongly acyclic and $E(\a)$ is therefore acyclic, forcing $m = 0$
and proving the claim.

We will be done if we can show 
that the negation of (i) implies the negation of (iii).
Let $E$ be a heap with property P2 but not property P1, and
suppose that $E$ is minimal with this property, in the sense that whenever
$\a$ is a maximal or minimal vertex of $E$ then $E(\a)$ has property P1.
(Note that $E(\a)$ will also have property P2 in this situation.)  
Any such heap $E$ is necessarily nontrivial, so by Lemma 2.3.1, $\Delta(E)$
also has property P2.

Let $\a$ be a minimal element of $E$.  The minimality property of $E$ shows
that $\ldescent{E} = \ldescent{E(\a)} \cup \{\a\}$, 
which implies that $\Delta(E(\a))
= \Delta(E) \backslash \{ \a \}$ (embedding $E$ in $\Delta(E)$ in the usual
way).  Furthermore, there is
no subheap $F$ of $E(\a)$ with $F \prec^+ E(\a)$ by Lemma 2.2.2 (ii), and 
$E(\a)$ has
property P1 by minimality of $E$.  Since any heap that has no balanced 
convex chains of length 2 and that is of the form $\a \circ T$ for $T$ trivial
must have property P1, $E(\a)$ is not a trivial heap.
By Lemma 2.2.3 (ii), $\Delta(E(\a)) = \d^m E_1$ (as elements of $TL(P, \C)$)
for some $m > 0$ and some trivial heap $E_1$.  This means that, as
elements of $TL(P, \C)$, 
$\Delta(E) \backslash \{ \a \} = \d^m E_1$ for some $m > 0$, which
contradicts (iii), completing the proof.
\qed\enddemo

\head 3. Classification of regular classes of heaps \endhead

In \S3, we shall prove Theorem 1.5.2.  The proof techniques we use
are reminiscent of those used in 
the classification of finite Coxeter groups (see
\cite{{\bf 11}, \S2}) and the classification of FC-finite Coxeter groups (see
\cite{{\bf 12}, \S4} or \cite{{\bf 5}, \S7}).

\subhead 3.1 Subgraphs and connected components \endsubhead

Our classification of regular classes of heaps is in terms of their concurrency
graphs (Definition 1.1.5).

\definition{Definition 3.1.1}
A graph $\Gamma$ is said to have {\it property R} if it is the 
concurrency graph of a regular class of heaps.
\enddefinition

Our aim is to classify all finite graphs with property R.  The key to the
procedure is the following observation.

\proclaim{Lemma 3.1.2}
Let $H(P, \C)$ be a class of heaps, let $P'$ be a subset of $P$ and let
$\C'$ be the restriction of $\C$ to $P'$.  
\item{\rm (i)}{There is a canonical inclusion,
$\iota$, of $H(P', \C')$ into $H(P, \C)$ that respects the partial order 
on heaps and the labelling function $\e : E \ra P$.}
\item{\rm (ii)}{The heap, $E$ has property P1 (respectively, property P2) if
and only if $\iota(E)$ does.}
\endproclaim

\demo{Proof}
Let $E = [E, \leq, \e]$ be a heap of $H(P', \C')$.  The heap $\iota(E)$ of
$H(P, \C)$ is $[E, \leq, \eta]$, where $(E, \leq)$ is the same poset as
before and $\eta$ is the obvious extension of $\e$ to $\C$.  If $\a,
\be$ are elements of $E$, then $\e(\a) \ \C' \ \e(\be)$ if and only if 
$\eta(\a) \ \C \ \eta(\be)$, so both axioms of Definition 1.1.1 hold, proving
(i).

The assertion of (ii) follows because properties P1 and P2 can be
defined using only the partial order on the heap and the function $\e$.
\qed\enddemo

\proclaim{Lemma 3.1.3}
If $\Gamma'$ is a full subgraph of a graph 
$\Gamma$ with property R, then $\Gamma'$ has property R.
\endproclaim

\demo{Proof}
If $\Gamma$ is the concurrency graph of $H(P, \C)$ and $\Gamma'$ is a full
subgraph of it, then $\Gamma'$ must correspond to $H(P', \C')$ for some
subset $P'$ of $P$, where $\C'$ is the restriction of $\C$ to $P'$.  If
we denote the canonical embedding $H(P', \C') \hookrightarrow H(P, \C)$ by
$\iota$, Lemma 3.1.2 (ii) 
shows that a heap $E = [E, \leq, \e]$ of $H(P', \C')$
has property P2 but not property P1 if and only if $\iota(E)$ does.  The
conclusion follows.
\qed\enddemo

Our considerations reduce quickly to problems about connected graphs, thanks
to the following result.

\proclaim{Lemma 3.1.4}
Let $H(P, \C)$ be a class of heaps with concurrency graph $\Gamma$, and suppose
that $\Gamma$ is the disjoint union of two subgraphs, 
$\Gamma_1$ and $\Gamma_2$.  Let $P_1$ and $P_2$ be the respective subsets of
$P$, and let $\C_1$ and $\C_2$ be the restrictions of $\C$ to $P_1$ and $P_2$
respectively.
\item{\rm (i)}{Any heap $E = [E, \leq, \e]$ of $H(P, \C)$ may be written
uniquely as the product $$
E = 
\iota_{P_1}(E_1) \circ \iota_{P_2}(E_2) = 
\iota_{P_2}(E_2) \circ \iota_{P_1}(E_1)
,$$ where $E_i$ is a heap in $H(P_i, \C_i)$ and $\iota_{P_i}(E_i)$ denotes the 
embedding of $E_i$ into $H(P, \C)$ as in Lemma 3.1.2 (i).}
\item{\rm (ii)}{The heap $E$ has property P1 (respectively, property P2) if
and only if $E_1$ and $E_2$ both have property P1 (respectively, property
P2).}
\endproclaim

\demo{Proof}
For $i \in \{1, 2\}$, define $E_i$ to be the heap of $H(P_i, \C_i)$ 
corresponding via the map $\iota_{P_i}$ to the subheap of $E$ whose vertices 
are precisely those vertices $\a$ with $\e(\a) \in P_i$.  Observe that if
$\a$ and $\be$ are vertices of $E$ with $\a \leq \be$ and 
$\e(\a) \ \C \ \e(\be)$, we must have $\a, \be \in P_i$ for the same $i$,
and $\e(\a) \ \C_i \ \e(\be)$.  It follows from Definition 1.1.1 that the
poset $(E, \leq)$ is the disjoint union of $(E_1, \leq_1)$ and $(E_2, \leq_2)$,
where $\leq_i$ is the restriction of $\leq$ to $E_i$.  

Claim (i) follows from the facts that 
(a) $k \circ l = l \circ k$ whenever it is not the case that
$\e(k) \ \C \ \e(l)$ and (b)
$E$ can be written as a finite $\circ$-product of singleton heaps.  For (ii),
we may replace the heaps $E_i$ by 
$\iota_{P_i}(E_i)$ by Lemma 3.1.2 (ii), and the claim follows by 
an argument similar to the last part of the proof of that result.
\qed\enddemo

\proclaim{Corollary 3.1.5}
A graph $\Gamma$ has property R if and only if all its connected components
do.
\endproclaim

\demo{Proof}
If $\Gamma$ has property R then any connected component does by Lemma 3.1.3.
If all the connected components have property R, then $\Gamma$ does by Lemma
3.1.4 (ii).
\qed\enddemo

\subhead 3.2 Some graphs with property R \endsubhead

In \S3.2, we prove that all the graphs mentioned in Theorem 1.5.2 have 
property R.  Much of the work for this has been done by Fan:

\proclaim{Theorem 3.2.1 (Fan)}
If $\Gamma$ is a graph of type $A_n$, $D_n$, $E_n$ or $\ti{A}_{n-1}$ with 
$n$ odd, then $\Gamma$ has property R.
\endproclaim

\demo{Proof}
For types $A$, $D$ and $E$, this is a restatement of 
\cite{{\bf 3}, Lemma 4.3.1},
and for type $\ti{A}_{n-1}$, this is a restatement of 
\cite{{\bf 4}, Proposition
3.1.2}.  (See also \cite{{\bf 9}, Theorem 3.4.1}.)
\qed\enddemo

\proclaim{Proposition 3.2.2}
If $\Gamma$ is a complete graph, $K_n$, then $\Gamma$ has property R.
\endproclaim

\demo{Proof}
Let $H(P, \C)$ be a class of heaps with concurrency graph $\Gamma$, and
let $E = [E, \leq, \e]$ be a heap of $H(P, \C)$ with property P2; we will
show that $E$ has property P1.

Let $T_1 \circ \cdots \circ T_p$ be the unique factorization of $E$ as
in Proposition 2.1.1.  Since $\Gamma$ is a complete graph, all the factors
$T_i$ must contain a single element, because the $T_i$ are trivial.  We
may assume $p > 1$ (or we are done).  Let $\a$ be the vertex in the factor
$T_1$, and let $\be$ be the vertex in the factor $T_2$.
Since $E$ has property P2, we cannot have $\e(\a) = \e(\be)$, but we must
have $\e(\a) \ \C \ \e(\be)$ because $\Gamma$ is complete.  We therefore
have $E(\a) \prec^- E$, and the claim follows by induction on $|E|$.
\qed\enddemo

We introduce the following definition for notational convenience.

\definition{Definition 3.2.3}
Let $E = [E, \leq, \e]$ be a heap of $H(P, \C)$ with unique factorization $$
E = T_1 \circ T_2 \circ \cdots \circ T_p
.$$  We say that the piece $w$ is {\it represented} in the factor $T_i$
by $\a \in E$ if there is a vertex $\a$ in the factor $T_i$ with $\e(\a) = w$.
We say that $w$ {\it occurs} in the factor $T_i$ if it is represented by some
vertex $\a$ in $T_i$.
\enddefinition

\example{Example 3.2.4}
In the heap of Example 1.1.2, the piece $3$ occurs in $T_1$ and $T_3$; it is 
represented by $b$ in $T_1$ and by $e$ in $T_3$.
\endexample

We now turn our attention to type $\ti{E}_6$.  It is notationally 
convenient to assign names to the vertices of the graph (\idest the pieces).  
We call the branch point $c$; its neighbours are $b$, $d$ and $f$,
and the other vertices adjacent to $b$, $d$ and $f$ are denoted
by $a$, $e$ and $g$, respectively.  (In Figure 1, the labels could read
$a$, $b$, $c$, $d$, $e$ along the top row, and then $f$ and $g$ reading
downwards.)  We will consider a minimal counterexample to regularity in
the sense of the proof of Theorem 1.5.1, that is, a heap $E$ with property
P2 but not P1, but such that any heap $E(\a)$ for a maximal or minimal vertex
$\a$ has property P1.  (By an argument similar to that used to prove Lemma
2.2.3 (iii), $E(a)$ will automatically have property P2.)

\proclaim{Lemma 3.2.5}
Let $H(P, \C)$ be a class of heaps with concurrency graph $\ti{E}_6$ and let
$E$ be a minimal heap with property P2 but not property P1 (in the sense of
the above discussion).  Let $T_1 \circ \cdots \circ T_p$ be the unique
factorization of $E$.  Then:
\item{\rm (i)}{$p \geq 3$;}
\item{\rm (ii)}{$c$ occurs in $T_2$;}
\item{\rm (iii)}{$b$, $d$ and $f$ are represented by vertices $\a_1, \a_2 
\in T_1$ and $\a_3 \in T_3$, in some order, and $\e(\a_3)$ does not occur in
$T_1$.}
\endproclaim

\demo{Proof}
Clearly $E$ cannot be trivial, because trivial heaps have property P1, so
$p > 1$.  If $E$ has property P2 but not P1, then the opposite heap $E^*$
does as well, and Lemma 2.2.1 (ii) makes it impossible to have $p = 2$ because
$\Gamma$ is a finite graph with no circuits.  This proves (i).

Choose a minimal vertex, $\gamma$, of $E$.  By the minimality of $E$, 
$E(\gamma)$
has property P1.  By Lemma 2.2.2 (ii) and repeated applications of 
Lemma 2.2.3 (iii) (starting with the heap $E(\gamma)$), we find that there
exists a heap $F$ with $F \prec^- E'$, where $$
E' = T_2 \circ T_3 \circ \cdots \circ T_p
.$$   Let $\alpha \in T_2$ and $\beta \in T_3$ be as in the definition of the
condition $F \prec^- E'$.  If $\kappa$ were to represent $\e(\beta)$ 
in $T_1$ then $\kappa < \alpha < \beta$ would be a balanced convex chain 
contradicting the assumption that $E$ has property
P2, so $\e(\beta)$ cannot occur in $T_1$.  However, Lemma 2.2.1 (ii) guarantees
that at least two neighbours of $\e(\alpha)$ in the graph $\Gamma$ occur in
$T_1$, so the valency of $\e(\alpha)$ must be at least 3, \idest 
$\e(\alpha) = c$,
proving (ii).  Part (iii) follows as a consequence of the same argument.
\qed\enddemo

\proclaim{Proposition 3.2.6}
A graph of type $\ti{E}_6$ has property R.
\endproclaim

\demo{Proof}
Let $H(P, \C)$ be a class of heaps with concurrency graph $\ti{E}_6$ and 
suppose that $E$ is a minimal heap with property P2 but not P1, as in the 
proof of Lemma 3.2.5.  Let $T_1 \circ \cdots \circ T_p$ be the unique
factorization of $E$.

By symmetry of $\Gamma$ and Lemma 3.2.5, we may assume that $b$ and $d$
occur in $T_1$ and $f$ occurs in $T_3$ but not $T_1$.  The only other vertex
that could occur in $T_1$ is $g$, since $T_1$ is trivial.  In $T_2$, there
cannot be any occurrences of $a$ or $e$, because then removal of the 
occurrences of $b$ or $d$ in $T_1$ would violate the hypotheses on $E$.
It follows that only $c$ can occur in $T_2$, and 
$|T_2| = 1$.  We have just established that the monomial of
$\text{\rm Co}(P, C)$ (where $C$ is the complementary relation of $\C$) 
corresponding to $T_1 \circ T_2$ is either $bdgc$ or $bdc$.

We know that $T_3$ contains an occurrence of $f$.  Any other elements of
$P$ occurring in $T_3$ must be adjacent to $c$ by Proposition 2.1.1.  However,
$b$ and $d$ cannot occur because subwords $bdgcb$ and $bdgcd$ are not allowed
in a heap with property P2, so $|T_3| = 1$.  Both the sequences given
correspond to heaps with property P1, so we must have $p > 3$.

Using a similar argument, we find that the only element of $P$ that can
occur in $T_4$ is $g$, and even this is not allowed if $g$ occurs in $T_1$,
because $E$ has property P2.  The monomial of $\text{\rm Co}(P, C)$ 
corresponding to $E$ can therefore be assumed to start $bdcfg\ldots$, and this
cannot be the complete monomial since the first 5 letters correspond to a
heap with property P1.  However, no further letters can be added to it on
the right without violating one of the hypotheses.  We conclude that no 
such heap $E$ exists, and that $\Gamma$ has property R.
\qed\enddemo

\subhead 3.3 Triangles, circuits and branch points \endsubhead

All that remains to prove Theorem 1.5.2 is to classify the connected graphs
that fail to have property R.  Lemma 3.1.3 reduces this problem to a problem
about complete subgraphs.  However, Lemma 3.1.3 is false if we delete 
the word ``full'' (see Proposition 3.2.2).  To get around this problem,
we first consider the case where $\Gamma$ contains a triangle.

\proclaim{Lemma 3.3.1}
If $\Gamma$ is a graph with property R, then $\Gamma$ has no full
subgraphs of either of the types shown in Figure 2.
\endproclaim

\topcaption{Figure 2} Some 4-vertex graphs without property R
\endcaption
\centerline{
\hbox to 3.500in{
\vbox to 1.250in{\vfill
        \includegraphics{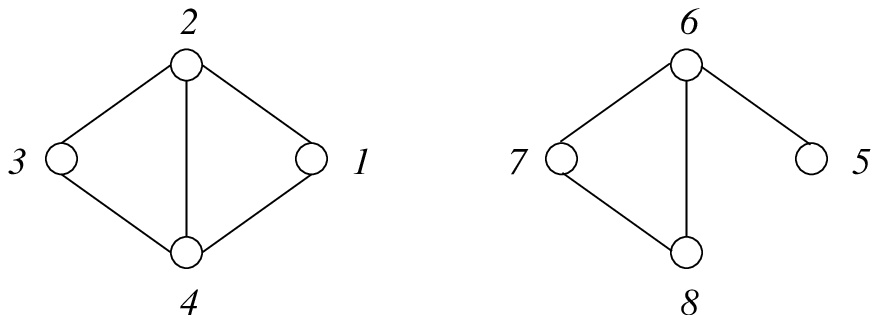}
}
\hfill}
}

\demo{Proof}
The two graphs shown in Figure 2 are the only connected graphs with four
vertices that contain a triangle but are not the complete graph $K_4$.

Consider the commutation monoid $\text{\rm Co}(P, C)$
corresponding to the graph on the left.
In this case, only the vertices 1 and 3 commute with each other, and the
monomial $(13)(2)(4)(13)$ corresponds to a heap with property P2 but not
property P1, so this graph does not have property R.  
More precisely, the corresponding heap is of the form
$T_1 \circ T_2 \circ T_3 \circ T_4$ where each factor corresponds to a
parenthetic expression in the monomial, so for example $T_3$ contains an
occurrence of the piece $4$, and no other vertices.

Now consider the graph on the right; this is a triangle on the vertices
6, 7 and 8, where 5 is connnected only to 6.  Then the monomial
$(58)(6)(7)(8)(6)(57)$ in the commutation monoid corresponds to a heap
(with the indicated unique factorization) that has property P2 but not
property P1.
\qed\enddemo

\remark{Remark 3.3.2}
A quick way to check that the monomials in the above proof have property
P2 is to verify that any two occurrences of the same generator $a$ are
separated by at least two other occurrences of generators not commuting with 
$a$ in $\text{\rm Co}(P, C)$.  Fan calls this property R3; see \cite{{\bf 3},
\S2}.
\endremark

\proclaim{Lemma 3.3.3}
Let $\Gamma$ be a finite connected graph containing a triangle, and suppose
that every 4-vertex connected full subgraph of $\Gamma$ that contains 
a triangle is the complete graph $K_4$.  Then $\Gamma$ is complete.
\endproclaim

\demo{Proof}
If $\Gamma$ has three vertices, there is nothing to prove, so we assume
$|\Gamma| > 3$.  Suppose $\Gamma$ is not complete, and let $\Gamma'$ be
a maximal complete subgraph of $\Gamma$, which must therefore be a proper
subgraph.  Let $\a$ be a vertex of $\Gamma \backslash \Gamma'$ 
adjacent to some vertex $\be$ of 
$\Gamma'$; this must exist since $\Gamma$ is connected.  Let $\g$ be
any other vertex of $\Gamma'$.  Since $|\Gamma'| \geq 3$, there is a triangle
in $\Gamma'$ containing $\be$, $\g$ and some other vertex, $\g'$.  Now
$\{\a, \be, \g, \g'\}$ induces a 4-vertex connected full subgraph of $\Gamma$
that contains a triangle, so it must be the complete graph.  In particular,
$\a$ is adjacent to $\g$, and, since $\g$ was arbitrary, the full subgraph
containing the vertices $\Gamma' \cup \{\a\}$ is a complete subgraph
properly containing $\Gamma'$, a contradiction.
\qed\enddemo

\proclaim{Corollary 3.3.4}
If $\Gamma$ is a finite connected graph with property R and $\Gamma$
contains a triangle (as a subgraph), then $\Gamma$ is complete.
\endproclaim

\demo{Proof}
Suppose $\Gamma$ is not complete.  Then by Lemma 3.3.3, $\Gamma$ contains
one of the graphs in Figure 2 as a full subgraph.  Lemma 3.3.1 and Lemma
3.1.3 then show that $\Gamma$ does not have property R.
\qed\enddemo

\proclaim{Lemma 3.3.5}
If $\Gamma$ is a finite incomplete connected graph with property R 
and $\Gamma$ contains a circuit, then $\Gamma$ is of type $\ti{A}_{n-1}$ 
for some odd number $n$.
\endproclaim

\demo{Proof}
Let $g_1, g_2, \ldots, g_{l-1}, g_l = g_1$ be a circuit of vertices in
$\Gamma$; note that $l > 3$.  Assume that $\Gamma$ is not itself a circuit.
We may assume without loss of
generality that $g_1$ is adjacent to some other vertex, $x$.
Since $\Gamma$ is incomplete, it contains no triangles by Corollary 3.3.4,
so $x$ is not adjacent either to $g_2$ or to $g_{l-1}$.  The monomial $$
(x g_{l-1}) (g_1) (g_2) \cdots (g_{l-2}) (g_{l-1}) (g_1) (x g_2)
$$ in the commutation monoid associated to $\Gamma$ then corresponds
to a heap with property P2 but not property P1.  This shows that $\Gamma$
is a circuit.

If $\Gamma$ is an even circuit, $g_1, g_2, \ldots, g_{2l}, g_{2l+1} = g_1$,
then the monomial $$
(g_1 g_3 g_5 \cdots g_{2l-1}) (g_2 g_4 g_6 \cdots g_{2l})
$$ corresponds to a heap with property P2 but not property P1.  This proves
that $n$ is odd.
\qed\enddemo

\proclaim{Lemma 3.3.6}
If $\Gamma$ is a finite incomplete connected graph with property R 
then $\Gamma$ contains at most one branch point.
\endproclaim

\demo{Proof}
Suppose $\Gamma$ contains two branch points, $c$ and $c'$.  
Then $\Gamma$ is not a circuit, and it cannot contain a circuit by 
Lemma 3.3.5.  Since $\Gamma$ is connected, it contains a shortest path $$
c = g_1, g_2, \ldots, g_{k-1}, g_k = c'
.$$  There are two vertices $x_1$ and $x_2$ distinct from $g_2$ that are
adjacent to $c$, and similarly there are two vertices $y_1$ and $y_2$ distinct
from $g_{k-1}$ that are adjacent to $c'$.  There are no coincidences among
the $g_i$, $x_i$ and $y_i$, because $\Gamma$ contains no circuits.  For the
same reason, $x_1$ is not adjacent to $x_2$, and $y_1$ is not adjacent to
$y_2$.  Then the commutation monoid associated to $\Gamma$ contains the
monomial $$
(x_1 x_2) (g_1) (g_2) \cdots (g_{k-1}) (g_k) (y_1 y_2)
,$$ which corresponds to a heap with property P2 but not property P1.
\qed\enddemo

\proclaim{Lemma 3.3.7}
If $\Gamma$ is a finite incomplete connected graph with property R 
then every vertex of $\Gamma$ has valency strictly less than $4$.
\endproclaim

\demo{Proof}
Suppose that $c$ is a vertex with valency $4$ or greater, and let $x_1, x_2,
x_3, x_4$ be four distinct vertices adjacent to $c$.  By Corollary 3.3.4, none
of the $x_i$ is adjacent to any other.  The monomial $$
(x_1 x_2) (c) (x_3 x_4)
$$ in the commutation monoid corresponds to a heap with property P2 but not
property P1.
\qed\enddemo

\subhead 3.4 The graphs $\Gamma(p, q, r)$ \endsubhead

\definition{Definition 3.4.1}
Let $p, q, r$ be nonnegative integers with $p \leq q \leq r$.  
We define $\Gamma(p, q, r)$ to be
the graph with $p + q + r + 1$ vertices containing one vertex, $c$, of 
valency 3 and disjoint arms of lengths $p$, $q$ and $r$ emanating from
$c$.
\enddefinition

\example{Example 3.4.2}
In Figure 1, the graph of type $A_n$ is $\Gamma(0, 0, n-1)$, the graph of
type $D_n$ is $\Gamma(1, 1, n-3)$ and the graph of type $E_n$ is
$\Gamma(1, 2, n-4)$.
\endexample

Lemmas 3.3.5, 3.3.6 and 3.3.7 have the following consequence.

\proclaim{Corollary 3.4.3}
If $\Gamma$ is a finite connected graph with property R then either
$\Gamma$ is complete, or $\Gamma$ is an $n$-gon for some odd $n$, or $\Gamma
= \Gamma(p, q, r)$ for some $p$, $q$ and $r$.
\qed\endproclaim

\proclaim{Lemma 3.4.4}
If $\Gamma(p, q, r)$ has property R, then $q < 3$.
\endproclaim

\topcaption{Figure 3} The graphs $\Gamma(1, 3, 3)$ and $\Gamma(2, 2, 3)$
\endcaption
\centerline{
\hbox to 3.138in{
\vbox to 2.597in{\vfill
        \includegraphics{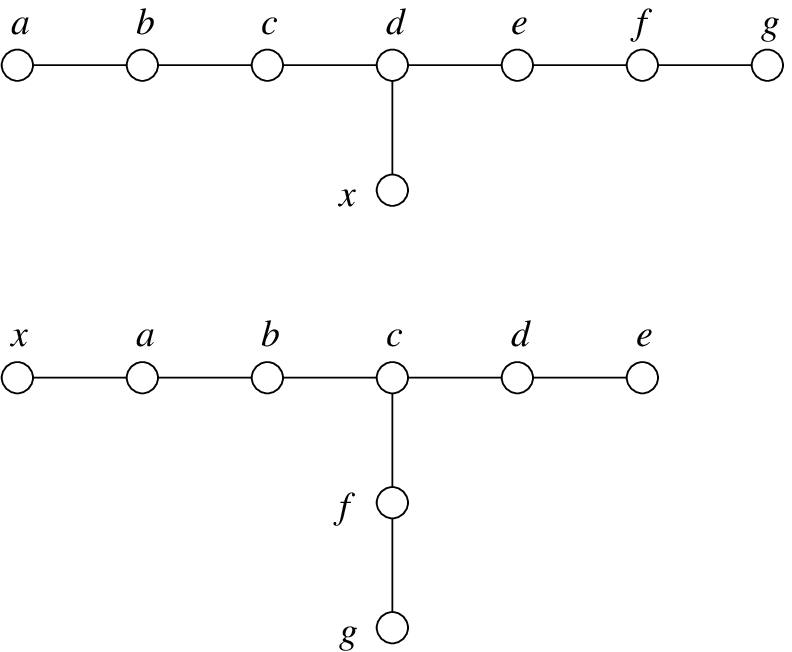}
}
\hfill}
}

\demo{Proof}
By Lemma 3.1.3, it is enough to show that $\Gamma(1, 3, 3)$ (also known
as $\ti{E}_7$) fails to have
property R.  Label the graph as shown in Figure 3.

The monomial $(acx)(bd)(ce)(df)(xeg)$ 
in the commutation monoid corresponds to a heap
with property P2 but not P1, and establishes the claim.
\qed\enddemo

\proclaim{Lemma 3.4.5}
If $\Gamma(p, q, r)$ has property R and $p = q = 2$, then $r = 2$.
\endproclaim

\demo{Proof}
By Lemma 3.1.3, it is enough to show that $\Gamma(2, 2, 3)$ fails to have
property R.  Label the graph as shown in Figure 3.

The monomial $(xbf)(ac)(bd)(ce)(df)(cg)(bf)(ac)(xbd)$ in the commutation
mon\-oid corresponds to a heap with property P2 but not P1, completing the
proof.
\qed\enddemo

\demo{Proof of Theorem 1.5.2}
We need only consider the case where $\Gamma$ is connected by Corollary 3.1.5.
In the light of Theorem 3.2.1 and Corollary 3.4.3, we only need to check 
that the claim
holds for $\Gamma = \Gamma(p, q, r)$.  If the hypotheses of Lemma 3.4.5 hold,
then we have type $\ti{E}_6$.  If $p = 0$, we are in type $A_n$.
Otherwise, if $q = 2$, we have $p = 1$, and
we are in type $E_n$.  By Lemma 3.4.4, the only other possibility is
$p = q = 1$, giving type $D_n$.
\qed\enddemo

\head 4. Concluding remarks \endhead

It might be interesting to investigate the representation theory of the
algebras $TL(P, \C)$ in the case where the concurrency graph has property R.
This was done by Fan in \cite{{\bf 3}} in the $ADE$ case, and there are
several papers on the case of type $\ti{A}_{n-1}$, including \cite{{\bf 7}} 
and \cite{{\bf 6}}.  More specifically, it would be interesting to know 
if the basis of monomials for $TL(P, \C)$ (Definition 1.4.7) is always a 
tabular basis in the sense of \cite{{\bf 8}}; this is true in type $ADE$ by
\cite{{\bf 8}, Theorem 4.3.5}, and in type $\ti{A}_{n-1}$ for $n$ odd by 
\cite{{\bf 8}, Theorem 6.4.8}.

\head Acknowledgements \endhead

The author thanks P.A. Haworth and A. Hulpke for helpful conversations.

\leftheadtext{}
\rightheadtext{}

\end